\documentclass[a4paper,12pt]{article}
\usepackage{amssymb,amsmath}
\newcommand{\meas}{{\rm meas}}

\newcommand{\R}{\mathbb{R}}

\newcommand{\ba}{\overline{a}}
\newcommand{\N}{\mathbb{N}}
\renewcommand{\b}[1]{{\bf #1}}
\renewcommand{\mod}[1]{\!\!\!\pmod{#1}}

\newcommand{\Z}{\mathbb{Z}}
\newcommand{\lb}[1]{\mbox{\scriptsize\bf #1}}
\newcommand{\twosum}[2]{\sum_{\substack{#1\\#2}}}
\newcommand{\nd}{\nmid}
\newcommand{\ep}{\varepsilon}
\newtheorem{theorem}{Theorem}
\newtheorem{lemma}{Lemma}
\newtheorem{conjecture}{Conjecture}

\begin{document}
\title{Pair Correlation for Fractional Parts of $\alpha n^2$}
\author{D.R. Heath-Brown\\Mathematical Institute, Oxford}
\date{}
\maketitle

\section{Introduction}

It was proved by Weyl \cite{Weyl} in 1916 that the sequence of values
of $\alpha n^2$ is uniformly distributed modulo 1, for any fixed 
real irrational $\alpha$.  Indeed Weyl's result covered sequences $\alpha n^d$
for any fixed positive integer exponent $d$.
However Weyl's work leaves open a number of questions concerning
the finer distribution of these sequences. It has been conjectured by
Rudnick, Sarnak and Zaharescu \cite{RSZ} that the fractional parts of
$\alpha n^2$ will have a Poisson distribution provided firstly that
$\alpha$ is
``Diophantine'', and secondly that if $a/q$ is any convergent to $\alpha$ then
the square-free part of $q$ is $q^{1+o(1)}$.  Here one says that
$\alpha$ is Diophantine if one has
\begin{equation}\label{Dio}
\left|\alpha-\frac{a}{q}\right|\gg_{\ep} q^{-2-\ep}
\end{equation}
for every rational number $a/q$ and any fixed $\ep>0$.  In particular
every real irrational algebraic number is Diophantine.  One would
predict that there are Diophantine numbers $\alpha$ for which the
sequence of convergents $p_n/q_n$ contains infinitely many squares
amongst the $q_n$.  If true, this would show
that the second condition is independent of the first.  Indeed one
would expect to find such $\alpha$ with bounded partial quotients.

The Poisson property can be phrased in terms of a sequence of
correlation functions.  We shall be concerned in the present paper
with the pair correlation function.  For a real sequence
$\theta=(\theta_n)_1^{\infty}$ considered modulo 1 
we define the pair correlation
function by
\[R(N,X;\theta):=N^{-1}\#\{m<n\le N:||\theta_m-\theta_n||\le
XN^{-1}\},\]
and if $\theta_n=\alpha n^2$ we write $R_{\alpha}(N,X)$ in place of
$R(N,X;\theta)$. 
If the sequence $\theta$ follows a Poisson distribution then we will
have
\begin{equation}\label{pc}
\lim_{N\rightarrow\infty} R(N,X;\theta)=X\;\;\;\mbox{for all}\; X>0.
\end{equation}
The statement (\ref{pc}) is in general weaker than the Poisson condition.
However we know rather little even about $R_{\alpha}(N,X)$.  For the
pair correlation it appears that one does not need the ``nearly 
square-free'' condition
for the numerators $q$ of the convergents $a/q$ of $\alpha$.  
We therefore make the following conjecture.

\begin{conjecture}
If $\alpha$ is Diophantine, then
\begin{equation}\label{pca}
R_{\alpha}(N,X)=N^{-1}\#\{m<n\le N:||\alpha(m^2-n^2)||\le
XN^{-1}\}\rightarrow X
\end{equation}
as $N\rightarrow \infty$, for each fixed $X\ge 0$.
\end{conjecture}

Related conjectures are already mentioned in the works of Rudnick and
Sarnak \cite{RS} and of Rudnick, Sarnak and Zaharescu \cite{RSZ}.

We remark that if $|\alpha-a/q|\le 1/(4q^3)$ infinitely often, then
(\ref{pca}) is false for every $X\in (1/4,1/2)$.  To see this one 
may take $N=q$
and note that the pairs $m$ and $n=q-m$ for $m<q/2$ will satisfy
$||\alpha(m^2-n^2)||\le 1/(4N)$, whence $R_{\alpha}(N,X)\ge 1/2+o(1)$.
Thus some condition on rational approximations to $\alpha$ will
clearly be necessary. A similar remark occurs in the paper of
Rudnick and Sarnak \cite{RS}.

Rudnick and Sarnak \cite{RS} were able to show that
``almost all'' $\alpha$, in the sense of Lebesgue measure, satisfy 
(\ref{pca}), but they remark that
they are not able to provide any explicit value of $\alpha$ which does
so.  An alternative proof of this result was given by Marklof and
Str\"{o}mbergsson \cite{MS}.
Our first result gives a third way to establish the ``almost-all'' property,
but more importantly it allows us to construct values of $\alpha$
for which (\ref{pca}) holds.

\begin{theorem}\label{con}
The statement {\rm (\ref{pca})} holds for almost all real $\alpha$.
Moreover there is a dense set of constructible values of $\alpha$ for
which {\rm (\ref{pca})} holds.
\end{theorem}

The second claim of the theorem deserves further comment.  What we
will do is to provide an informal algorithm, which, for any closed
interval $I$ of positive length, provides a convergent 
sequence of rational numbers
belonging to $I$, whose limit $\alpha$ satisfies (\ref{pca}). It could
be said Rudnick and Sarnak were hoping for a more explicit
construction, akin to that for Liouville numbers, for example.
However, from a logical point of view there is no difference
between our construction and that of other more familiar real
numbers.  The reader might also feel happier if we had given an
explicit example of an admissible $\alpha$, by displaying its decimal
expansion; but since our construction provides a dense set of values,
that would be uninformative.  We can safely assert that
\[\alpha=3.14159265358\ldots\]
satisfies (\ref{pca}), but this will not help the reader's intuition!

We note at this point that our proof of Theorem \ref{con} provides
slightly more. Indeed there is a positive constant $\eta$ (we may take
$\eta=1/200$) such that for almost all $\alpha$, and in particular for
those $\alpha$ which we construct, we have
$R_{\alpha}(N,X)=X+O(N^{-\eta})$ uniformly for $N^{-\eta}\le X\le
N^{\eta}$.

Rudnick and Sarnak \cite{RS} proved that (\ref{pc}) holds
for the sequence $(\alpha n^d)_1^{\infty}$ for almost all $\alpha$, for every
$d\ge 2$.  However our approach appears to work only for $d=2$.

In proving Theorem \ref{con} we shall show that (\ref{pca}) holds for
all $\alpha$ satisfying three conditions, which are explained in
detail in \S 3.  The first of these is that
\[\left|\alpha-\frac{a}{q}\right|\gg q^{-2-1/200}\]
for every approximation $a/q$ to $\alpha$.  It is of interest 
that this requirement is not quite as strong as (\ref{Dio}).  The
second condition is roughly that if $a_n/q_n$ are the continued
fraction convergents to $\alpha$, then $q_n$ is ``almost odd and
square-free''.  It seems conceivable that one could adapt the proof to
avoid this condition.  The third assumption on $\alpha$ is that $a_n$
does not lie in a certain small ``bad'' set $B(q_n)$, if $n$ is large
enough. One would conjecture that the sets $B(q)$ are empty for all
sufficiently large $q$.  Thus in this approach it is the sets $B(q)$
which are the real stumbling block in any attack on Conjecture 1.

A related approach to Conjecture 1 has been investigated by Truelsen
\cite{Tr}.  This is based on a hypothesis concerning the average
value of the function
\[\tau^*_M(n):=\#\{(a,b)\in\N^2:a,b\le M, ab=n\}\]
in short arithmetic progressions.  Such a hypothesis is related to the
condition giving our bad sets $B(q)$.  Truelsen proves that his
hypothesis holds on average, in a suitable sense.  Our Lemma 3 is in a
similar vein, but the two results are not directly comparable.
\bigskip

Our second result gives partial support to Conjecture 1, by
describing the behaviour of $R_{\alpha}(N,X)$ as $X$ grows.

\begin{theorem}\label{main}
Suppose that $\alpha\in\R$ and $\kappa>1$ satisfy
\[\left|\alpha-\frac{a}{q}\right|\ge \frac{1}{\kappa q^{9/4}}\]
for all fractions $a/q$.  Then
\[R_{\alpha}(N,X)=X+O(X^{7/8})+O(\kappa^2(\log N)^{-1})\]
uniformly in all the parameters, for $1\le X\le \log N$.
\end{theorem}

This result applies in particular whenever $\alpha$ is Diophantine.
It shows that, in the limit as $N\rightarrow\infty$, the function 
$R_{\alpha}(N,X)$ is approximately equal to $X$ for large
$X$. Moreover we have the correct order of magnitude
\[X\ll_{\kappa} R_{\alpha}(N,X)\ll_{\kappa} X\]
as soon as $X\gg_{\kappa} 1$.  

For a fixed $X\le 1$ we are unable to prove even that
$R_{\alpha}(N,X)\gg 1$ in general.  However the method 
used to establish Theorem
\ref{main} can be adapted to yield some non-trivial upper bounds, of
the form $R_{\alpha}(N,X)\ll X^{\theta}$ with $\theta>0$.  Here we require
$\alpha$ to be Diophantine, and $(\log N)^{-\delta}\le X\le 1$ for a
suitably small constant $\delta>0$.

In discussing Theorem \ref{main} it is natural to
examine the case $\alpha=a/q$, which leads to consideration of
congruences $a(m^2-n^2)\equiv r\mod{q}$ with $r$ small.  Thus it would
be interesting to know about the number of solutions $u,v\le N$ of
$uv\equiv c\mod{q}$, for a fixed $c$.  During the proof of Theorem
\ref{main} we will use a result of Linnik
and Vinogradov \cite{LV} which shows that
\[\twosum{k\le N^2}{k\equiv c\mod{q}}d(k)\ll_{\delta}
\phi(q)q^{-2}N^2\log N\]
uniformly for $(c,q)=1$ and $q\le N^{2(1-\delta)}$, for any fixed
$\delta\in(0,1)$.  (Here $d(k)$ is the divisor function.)
However for our problem we expect that the factor $\log N$ can be
removed, and we make the following conjecture.
\begin{conjecture}
For any fixed $\delta\in(0,1)$ we have
\[\#\{(u,v)\in\N^2:\, u,v\le N,\, uv\equiv c\mod{q}\}\ll_{\delta}
\phi(q)q^{-2}N^2\]
uniformly for $(c,q)=1$ and $q\le N^{2(1-\delta)}$.
\end{conjecture}
Unfortunately it appears that the techniques used by Linnik and
Vinogradov do not work for the above variant of their problem.  It is
no coincidence that Conjecture 2 can be reformulated using Truelsen's
function $\tau^*_N(n)$.

In order to put Theorem \ref{main} into context it may be helpful to
record what one can say about arbitrary sequences
$\theta$.  For this purpose it will be more convenient
to use a weighted pair correlation function
\[R_0(N,X;\theta):=N^{-1}\sum_{m,n\le N}
\left\{1-\frac{||\theta_m-\theta_n||}{X/N}\right\}^{+}.\]
\begin{theorem}\label{triv}
Let $\theta=(\theta_n)_1^{\infty}$ be an arbitrary real sequence.  
\begin{enumerate}
\item[(i)] We have {\rm (\ref{pc})} if and only if
\begin{equation}\label{pc0}
\lim_{N\rightarrow\infty} R_0(N,X;\theta)=1+X\;\;\;
\mbox{for all}\; X>0.
\end{equation}
\item[(ii)] We have $R_0(N,X;\theta)\ge \max(1,X)$ for all $X\ge 0$.
\item[(iii)] We have $R_0(N,X+Y;\theta)\le R_0(N,X;\theta)+
R_0(N,Y;\theta)$ for all $X,Y\ge 0$.
\end{enumerate}
\end{theorem}
Notice in particular that in part (i) we make no assumption about
uniformity with respect to $X$ in either of the limits involved.
We remark that Part (ii) can be strengthened slightly with a little
more work.  If
$X=[X]+\xi$, where $[X]$ is the integer part of $X$, then
\[R_0(N,X;\theta)\ge X+\frac{\xi-\xi^2}{X}\]
for $X>0$.  Moreover we have equality whenever the sequence $\theta$
consists of equally spaced points.  Notice here that 
$X+(\xi-\xi^2)/X\ge\max(1,X)$ for $X>0$.

Part (ii) shows that $R(N,X;\theta)\ge X+O(1)$ on average with respect
to $X$ (by virtue of (\ref{add})).  Thus the lower bound
implicit in Theorem \ref{main} holds, on average, for any sequence $\theta$.
Moreover, we have
\[\frac{R_0(N,X;\theta)-1}{2}\le R(N,X;\theta)\le 
R_0(N,2X;\theta)-1.\]
Hence Theorem \ref{triv} 
shows that for any sequence $(\theta_n)_1^{\infty}$ one has 
\[X\ll R(N,X;\theta)\ll X\]
for $X\ge 2$, say, providing only that
$R(N,1;\theta)\ll 1$.
\bigskip

The author was introduced to the topic of this article by Jimi Truelsen.  
His input, through a number of interesting conversations, and his
helpful comments on an earlier draft of this paper,
is gratefully acknowledged.

\section{Proof of Theorem \ref{triv}}

In this section we give the rather easy proof of Theorem \ref{triv}.
For part (i) we first show that (\ref{pc}) implies (\ref{pc0}).  We 
use the fact that
\begin{equation}\label{add}
R_0(N,X;\theta)=1+\frac{2}{X}\int_0^X R(N,t;\theta)dt.
\end{equation}
For any fixed positive integer $K$ we have
\[\int_0^X R(N,t;\theta)dt\le 
\frac{X}{K}\sum_{k=1}^{K}R(N,\frac{kX}{K};\theta),\]
since $R(N,t;\theta)$ is non-decreasing with respect to $t$.
We now let $N$ tend to infinity and apply (\ref{pc}), whence
\[\limsup_{N\rightarrow\infty}\int_0^X
R(N,t;\theta)dt\le
\frac{X}{K}\sum_{k=1}^{K}\frac{kX}{K}=\frac{X^2(K+1)}{2K}\]
for any fixed positive integer $K$.  Since $K$ is arbitrary it follows
that
\[\limsup_{N\rightarrow\infty}\int_0^X
R(N,t;\theta)dt\le\frac{X^2}{2}.\]
The corresponding lower bound
\[\liminf_{N\rightarrow\infty}\int_0^X
R(N,t;\theta)dt\ge\frac{X^2}{2}\]
follows similarly from the inequality
\[\int_0^X R(N,t;\theta)dt\ge 
\frac{X}{K}\sum_{k=0}^{K-1}R(N,\frac{kX}{K};\theta).\]
Thus (\ref{pc}) implies (\ref{pc0}).  

To establish the reverse implication we use a standard Tauberian
argument.  For any $\Delta>0$ we have
\begin{eqnarray*}
R(N,X;\theta)&\le&\Delta^{-1}\int_X^{X+\Delta}R(N,t;\theta)dt\\
&=&\frac{X+\Delta}{2\Delta}\{R_0(N,X+\Delta;\theta)-1\}-
\frac{X}{2\Delta}\{R_0(N,X;\theta)-1\}.
\end{eqnarray*}
We let $N$ tend to infinity and apply (\ref{pc0}) to obtain
\[\limsup_{N\rightarrow\infty}R(N,X;\theta)\le
\frac{X+\Delta}{2\Delta}(X+\Delta)-\frac{X}{2\Delta}X=
X+\frac{\Delta}{2}.\]
Since $\Delta>0$ was arbitrary we deduce that
\[\limsup_{N\rightarrow\infty}R(N,X;\theta)\le X.\]
The corresponding lower bound is trivial if $X=0$ and otherwise
follows as above, starting with the fact that
\begin{eqnarray*}
R(N,X;\theta)&\ge&\Delta^{-1}\int_{X-\Delta}^{X}R(N,t;\theta)dt\\
&=&\frac{X}{2\Delta}\{R_0(N,X;\theta)-1\}-
\frac{X-\Delta}{2\Delta}\{R_0(N,X-\Delta;\theta)-1\}
\end{eqnarray*}
for any $\Delta\in(0,X)$.  This establishes part (i) of Theorem
\ref{triv}.

For the second part of the theorem we use the fact that
\begin{equation}\label{G}
XR_0(N,X;\theta)=\int_0^1 L(t,X)^2dt
\end{equation}
where
\begin{equation}\label{GG}
L(t,X)=\#\{n:||\theta_n-t||\le (2N)^{-1}X\}.
\end{equation}
However
\begin{eqnarray*}
\int_0^1 L(t,X)dt&=&
\int_0^1\#\{n:||\theta_n-t||\le (2N)^{-1}X\}\,dt\\
&=&\sum_{n=1}^N\,\int_{\theta_n-(2N)^{-1}X}^{\theta_n+(2N)^{-1}X}dt\\
&=&\sum_{n=1}^N N^{-1}X\\
&=&X,
\end{eqnarray*}
whence Cauchy's inequality yields
\[X^2\le\int_0^1 L(t,X)^2dt=XR_0(N,X;\theta),\]
from which assertion (ii) of Theorem \ref{triv} follows, since the
terms $m=n$ yield $R_0(N,X;\theta)\ge 1$.

To handle part (iii) we again use (\ref{G}). We have
\begin{eqnarray*}
L(t,X+Y)&\le& L(t+(2N)^{-1}Y,X)+L(t-(2N)^{-1}X,Y)\rule[-3mm]{0mm}{1mm}\\
&=&\sqrt{X}\, \frac{L(t+(2N)^{-1}Y,X)}{\sqrt{X}}+
\sqrt{Y}\, \frac{L(t-(2N)^{-1}X,Y)}{\sqrt{Y}},
\end{eqnarray*}
so that
\[L(t,X+Y)^2\le
(X+Y)\left\{\frac{L(t+(2N)^{-1}Y,X)^2}{X}+
\frac{L(t-(2N)^{-1}X,Y)^2}{Y}\right\}\]
by Cauchy's inequality.  If we apply (\ref{G}) to each term this
produces
\[(X+Y)R_0(N,X+Y;\theta)\le
(X+Y)\{X\frac{R_0(N,X;\theta)}{X}+Y\frac{R_0(N,Y;\theta)}{Y}\},\]
which gives the required inequality.

\section{Proof of Theorem \ref{con}}

In this section we begin by presenting a proof that (\ref{pca}) 
holds for almost all real $\alpha$.  This is an immediate consequence 
of the following result, via the Borel--Cantelli Lemma.  Of course it
suffices to consider $\alpha\in[0,1]$, since the set of $\alpha$ for
which (\ref{pca}) holds has period 1.
\begin{lemma}\label{con1}
There is an explicit sequence of open intervals $I_n$, with
\[\sum_{n=1}^{\infty}\meas(I_n)<\infty,\]
such that if $\alpha\in [0,1]$ lies in only finitely many of
the $I_n$ then {\rm (\ref{pca})} holds for each fixed $X$.
\end{lemma}

The proof of this result will take up the bulk of this section.  It
will become clear in the course of this work, just what is meant by
the word ``explicit'' in the statement above.  At the end of this
section we shall show how Lemma \ref{con1} allows us to construct
values of $\alpha$ for which (\ref{pca}) holds.

In the course of the proof
we shall use a small parameter $\eta>0$, which we shall take to be 
\[\eta=1/200\]
However we prefer to use the notation $\eta$, which will make it
clearer why it suffices to use any sufficiently small positive value.
The reader will easily confirm at each step that $\eta=1/200$ is indeed
suitably small.

The intervals $I_n$ which we produce will be of three types.  We begin
by including all intervals
\[\left(\frac{a}{q}-\frac{1}{[q^{2+\eta}]}\, , \,
\frac{a}{q}+\frac{1}{[q^{2+\eta}]}  \right),\;\;\; 0\le a\le q\]
 among the $I_n$.  Clearly, for these we have
\[\sum_n \meas(I_n)\ll \sum_{q=1}^{\infty}\sum_{a=0}^q q^{-2-\eta}\ll
\sum_{q=1}^{\infty} q^{-1-\eta}<\infty.\]
If $\alpha$ belongs to only finitely many of these intervals we will have
\[\left|\alpha-\frac{a}{q}\right|\ge q^{-2-\eta}\]
for $q\ge q_0(\alpha)$, say.  Then if $a_j/q_j$ and
$a_{j+1}/q_{j+1}$ are successive convergents to $\alpha$, with $q_j>
q_0$, we will have
\[\frac{1}{q_j q_{j+1}}\ge\left|\alpha-\frac{a_j}{q_j}\right|\ge
q_j^{-2-\eta},\]
whence $q_{j+1}\le q_j^{1+\eta}$.  It follows that if $N\ge q_0$ there will
be a convergent $a/q$ with
\begin{equation}\label{Nb}
N^{3/2+\eta}\le q\le N^{3/2+3\eta}
\end{equation}
and hence
\[q^{2/3-4\eta/3}\le N\le q^{2/3-\eta/3}. \]

Since $a/q$ is a convergent to $\alpha$ we will have
\[\left|\alpha-\frac{a}{q}\right|\le q^{-2}.\]
Let $\alpha=a/q+\phi$, so that $|\phi|\le q^{-2}$.  We will find it
convenient to write $\ba$ for the multiplicative inverse of 
$a$ modulo $q$.  We begin our analysis of $R_{\alpha}(N,X)$
by observing that
\[2R_{\alpha}(N,X)+1=N^{-1}\{m,n\le N:||\alpha(m^2-n^2)||\le X/N\}.\]
Now, if $m^2-n^2\equiv \ba r\mod{q}$ with $|r|\le
Xq/N-N^2/q$, then
\[||\alpha(m^2-n^2)||\le ||\frac{a}{q}(m^2-n^2)||+|\phi(m^2-n^2)|
\le |r|q^{-1}+q^{-2}N^2\le XN^{-1}.\]
Similarly, if $||\alpha(m^2-n^2)||\le X/N$ then $m^2-n^2\equiv \ba 
r\mod{q}$ with 
\[|r|\le Xq/N+N^2/q.\]
We shall write
\[A(N,q,c)=\#\{m,n\le N: m^2-n^2\equiv c\mod{q}\},\]
whence
\begin{eqnarray*}
N^{-1}\sum_{|r|\le Xq/N-N^2/q}A(N,q,\ba r)&\le&
2R_{\alpha}(N,X)+1\\
&\le &N^{-1}\sum_{|r|\le Xq/N+N^2/q}A(N,q,\ba r).
\end{eqnarray*}
If we now impose the condition that
\[N^{-\eta}\le X\le N^{\eta}\]
then $0<Xq/N-N^2/q<Xq/N+N^2/q<q$ for large enough $N$, in view of (\ref{Nb}).
Moreover, if $q\nd c$ we have
\[A(N,q,c)\ll \twosum{|k|\le N^2}{k\equiv c\mod{q}}d(|k|)\le
N^{2+\eta}q^{-1},\]
since the case $k=0$ cannot occur.  It follows that
\[2R_{\alpha}(N,X)+1=
N^{-1}\sum_{|r|\le Xq/N}A(N,q,\ba r)+O(N^{3+\eta}q^{-2}).\]
The error term here is $O(N^{-\eta})$, and
\[A(N,q,0)=N+O(\sum_{k\le N^2, q|k}d(k))=N+O(N^{2+\eta}q^{-1}),\]
whence
\begin{equation}\label{R}
R_{\alpha}(N,X)=
N^{-1}\sum_{1\le r\le Xq/N}A(N,q,\ba r)+O(N^{-\eta}).
\end{equation}

We shall see that the expected value of $A(N,q,c)$ is about
$(N/q)^2A_0(q,c)$, where
\[A_0(q,c):=\#\{1\le m,n\le q:\,m^2-n^2\equiv c\mod{q}\}.\]
We write
\[\Delta(M,q,c):=|A(M,q,c)-(\frac{M}{q})^2A_0(q,c)|\]
and
\[\Delta^*(q,c):=\max_{M\le q^{2/3}}\Delta(M,q,c).\]
With this notation the technical result which is the key to our 
approach is the following.
\begin{lemma}\label{disp}
Let $q$ be a positive integer.  Write $q$ as a product of prime powers
in the form $q=\prod_p p^{e(p)}$ and set
\[q_1:=\prod_{p=2\mbox{ \scriptsize{or} } e(p)>1}p^{e(p)}.\]
Then
\[\sum_{c=1}^q\Delta^*(q,c)^2
\ll q^{3/2+4\eta}q_1^{3},\]
with an implied constant which is effectively computable in terms of
$\eta$.
\end{lemma}
We will prove this later, in \S \ref{pdisp}.

In view of Lemma \ref{disp} we will include among the intervals $I_n$
described in Lemma \ref{con1} a second category, namely all those 
\[I_n=\left(\frac{a}{q}-\frac{1}{[q^{2}]}\, , \,
\frac{a}{q}-\frac{1}{[q^{2}]}  \right),\;\;\;(0\le a\le q)\]
for which $q_1\ge q^{2\eta}$.  Since 
\[1\le \frac{q_1^{1/2}}{q^{\eta}}\]
whenever $q_1\ge q^{2\eta}$ we see that for these intervals we have
\[\sum_n \meas(I_n)\ll \sum_{q:\, q_1\ge q^{2\eta}}\;\sum_{a=1}^q q^{-2}
\ll \sum_{q=1}^{\infty}\frac{q_1^{1/2}}{q^{\eta}}q^{-1}\]
\[=\{1+\frac{2^{1/2}}{2^{1+\eta}}
+\frac{2}{4^{1+\eta}}+\frac{2^{3/2}}{8^{1+\eta}}+\ldots\}
\prod_{p>2}\{1+\frac{1}{p^{1+\eta}}+
\sum_{e=2}^{\infty}\frac{p^{e/2}}{p^{e(1+\eta)}}\}< \infty.\]

From now on we may assume that we have $q_1\le q^{2\eta}$ for all
values of $q$ under consideration, so that the estimate in Lemma
\ref{disp} is of order $q^{3/2+10\eta}$.
We now define a set $B(q)$ of ``bad'' values for $a$ by setting
\[B(q):=\{0\le a\le q:\,\sum_{r\le q^{1/3+2\eta}}\Delta^*(q,\ba r)\ge 
q^{2/3-2\eta}\}.\]
Then Lemma \ref{disp} yields
\begin{eqnarray*}
q^{2/3-2\eta}\#B(q)&\ll& \sum_{a=0}^q\sum_{r\le
  q^{1/3+\eta}}\Delta^*(q,\ba r)\\
&\ll& q^{1/3+2\eta}\sum_{c=1}^q\Delta^*(q,c)\\
&\ll& q^{1/3+2\eta}.q^{1/2}.(q^{3/2+10\eta})^{1/2},
\end{eqnarray*}
by Cauchy's inequality.  Thus
\[\#B(q)\ll q^{11/12+9\eta}.\]
To handle the bad values of $a$ we introduce our third class of
intervals $I_n$, defined as
\[I_n=\left(\frac{a}{q}-\frac{1}{[q^2]}\, , \,
\frac{a}{q}+\frac{1}{[q^2]}  \right),\;\;\; a\in B(q),\]
and observe that
\begin{eqnarray*}
\sum_n \meas(I_n)&\ll& 
\sum_{q=1}^{\infty}\;\sum_{0\le a\le q,\, a\in B(q)}q^{-2}\\
&\ll& \sum_{q=1}^{\infty}q^{-13/12+9\eta}<\infty,
\end{eqnarray*}
providing that we choose $\eta<1/108$.  Thus for the three classes of
intervals we have defined we have
\[\sum_n \meas(I_n)<\infty\]
providing that we choose $\eta=1/200$, say.

When $a\not\in B(q)$ the estimate (\ref{R}) produces
\begin{eqnarray*}
R_{\alpha}(N,X)&=&Nq^{-2}\sum_{r\le Xq/N}A_0(q,\ba r)\\
&&\hspace{2cm}\mbox{}
+O(N^{-1}\sum_{r\le q^{1/3+2\eta}}\Delta^*(q,\ba r))+O(N^{-\eta})\\
&=&Nq^{-2}\sum_{r\le Xq/N}A_0(q,\ba r)+O(N^{-1}q^{2/3-2\eta})
+O(N^{-\eta})
\end{eqnarray*}
for $q^{2/3-4\eta/3}\le N\le q^{2/3-\eta/3}$.  Under this condition the
two error terms are both $O(N^{-\eta})$.

We proceed to investigate
\[\sum_{r\le R}A_0(q,\ba r).\]
The function $A_0(q,r)$ is multiplicative with respect to $q$.  Thus
if $q_1$ is defined as in Lemma \ref{disp} and $q_0:=q/q_1$ we will
have $A_0(q,r)=A_0(q_0,r)A_0(q_1,r)$.  Since $q_0$ is odd we have
\[A_0(q_0,r)=\#\{u,v\le q_0:\, uv\equiv r\mod{q_0}\},\]
whence $A_0(q_0,kr)=A_0(q_0,r)$ whenever $k$ is coprime to $q_0$.
Thus
\[\sum_{r\le R}A_0(q,\ba r)=
\sum_{d|q_0}\sum_{s=1}^{q_1}A_0(q_0,d)A_0(q_1,s)U(R,q_0,q_1:d,s)\]
where
\begin{eqnarray*}
U(R,q_0,q_1:d,s)&=&
\#\{r\le R:\, (q_0,r)=d,\,\ba r\equiv s\mod{q_1}\}\\
&=&\sum_{e|q_0/d}\mu(e)\#\{r\le R:\, de|r,\,\ba r\equiv s\mod{q_1}\}\\
&=&\sum_{e|q_0/d}\mu(e)\{\frac{R}{deq_1}+O(1)\}.\\
\end{eqnarray*}
It follows that
\[\sum_{r\le R}A_0(q,\ba r)=R\Sigma+O(E),\]
say, where
\[\Sigma=\sum_{d|q_0}\sum_{e|q_0/d}\mu(e)
\sum_{s=1}^{q_1}A_0(q_0,d)A_0(q_1,s)/deq_1\]
and
\[E=\sum_{d|q_0}\sum_{e|q_0/d}\sum_{s=1}^{q_1}A_0(q_0,d)A_0(q_1,s).\]
In a precisely similar way we have
\[\sum_{r=1}^{q}A_0(q,r)=q\Sigma+O(E),\]
with the same values of $\Sigma$ and $E$, so that
\[\sum_{r\le R}A_0(q,\ba r)=\frac{R}{q}\sum_{r=1}^{q}A_0(q,r)+O(E)\]
for $R\le q$.

Trivially we have
\[\sum_{s=1}^{q_1}A_0(q_1,s)=q_1^2.\]
Moreover the congruence $x^2\equiv k\mod{q_0}$ has $O(q_0^{\eta})$ 
solutions $x\mod{q_0}$, uniformly in $k$, since $q_0$ is square-free.
Thus $A_0(q_0,d)\ll q_0^{1+\eta}$.  These bounds show that
$E\ll q_0^{1+2\eta}q_1^2\ll q^{1+2\eta}q_1\ll q^{1+4\eta}$.  Finally
\[\sum_{r=1}^{q}A_0(q,r)=q^2,\]
whence
\[\sum_{r\le R}A_0(q,\ba r)=Rq+O(q^{1+4\eta})\]
for $R\le q$.  We therefore deduce that
\[R_{\alpha}(N,X)=X+O(N^{-\eta})+O(Nq^{-1+4\eta}).\]
The final error term will be $O(N^{-\eta})$ for 
$q^{2/3-4\eta/3}\le N\le q^{2/3-\eta/3}$.

In conclusion we have shown that (\ref{pca}) holds uniformly for
$N^{-\eta}\le X\le N^{\eta}$, providing that $\alpha$ lies in none of
the intervals $I_n$.  This establishes Lemma \ref{con1}.
\bigskip

To complete the proof of Theorem \ref{con} we describe an algorithm
which will generate explicit values of $\alpha$ for which (\ref{pca})
holds.  Suppose we are given a closed interval $I$ of positive length,
and that we wish to construct a suitable $\alpha$ belonging to $I$.
Without loss of generality we may assume that $I\subseteq [0,1]$,
since the property (\ref{pca}) has period 1 in $\alpha$.  Moreover we
shall assume that $I$ has rational end-points, as we clearly
may. Finally we write $L$ for
the length of $I$.  Now consider the following algorithm.

We begin by taking $\eta=1/200$, and we compute an integer $N$ such that
\[\sum_{n=N}^{\infty}I_n<L/2.\]
We have not specified a numbering for the intervals $I_n$, but it
would be easy to do so. The contribution from the first two classes 
of intervals $I_n$ is relatively easy to calculate.  For the third class
one would need to make explicit the implied constant in Lemma
\ref{disp}, but there is no theoretical difficulty in doing this.

Now, for each integer $k>N$ define
\[F_k:=I\setminus\bigcup_{n=N}^k I_n.\]
This is a finite union of closed intervals with rational end points,
since the intervals $I_n$ also had rational end points.
It is important to notice here that $F_k$ cannot be empty, since
\[\meas(F_k)\ge \meas(I)-\sum_{n=N}^k\meas(I_n)>L/2>0.\]
We may then compute the set of end points of all the intervals which make up
$F_k$, and take $r_k$
to be the smallest such end point.  Thus $r_k$ is an explicitly
computable rational number, with
\[r_k\in F_k\subseteq I.\]

It is clear from the definition that the sets $F_k$ are nested, with
$F_N\supseteq F_{N+1} \supseteq F_{N+2}\ldots$, whence the
sequence $r_k$ must be non-decreasing. It follows that it converges to
a limit, $\alpha$ say.  Take any integer $j\ge N$.  Then, 
since $r_k\in F_k\subseteq F_j$ for all
$k\ge j$, and $F_j$ is closed, it follows that $\alpha\in F_j$.
However this holds for all $j\ge N$, whence 
\[\alpha\in\bigcap_{j=N}^{\infty}F_j=I\setminus\bigcup_{n=N}^{\infty} I_n.\]
We therefore see that $\alpha$ lies in none of the $I_k$ for $k\ge N$,
so that (\ref{pca}) holds for $\alpha$, for all $X$.

The $\alpha\in I$ that we have produced has been ``constructed'' in the
sense that we have given a procedure for determining a sequence of
rationals which converges to $\alpha$.  This completes the proof of
Theorem \ref{con}

\section{Proof of Lemma \ref{disp}}\label{pdisp}

We begin this section by considering
\[S:=\sum_{c=1}^q \Delta(N,q,c)^2,\]
for which we prove the following result.
\begin{lemma}\label{dispa}
Let $q$ be a positive integer, and let $q_1$ be defined as in Lemma
{\rm \ref{disp}}.  Then if $N\le q^{2/3}$ we have
\[\sum_{c=1}^q \Delta(N,q,c)^2\ll_{\eta}q^{4/3+4\eta}q_1^{3},\]
with an implied constant which is effectively computable in terms of
$\eta$.
\end{lemma}

We start by observing that
\begin{eqnarray}\label{kk}
S&=&\sum_c A(N,q,c)^2-2N^2q^{-2}\sum_c A(N,q,c)A_0(q,c)+
N^4q^{-4}\sum_c A_0(q,c)^2\nonumber\\
&=&S_1-2N^2q^{-2}S_2+N^4q^{-4}S_3,
\end{eqnarray}
say.  Clearly
\[S_1=\#\{x_1,x_2,x_3,x_4\le N:\,q|x_1^2+x_2^2-x_3^2-x_4^2\},\]
\[S_2=\#\{x_1,x_3\le N,\; x_2,x_4\le q:\,q|x_1^2+x_2^2-x_3^2-x_4^2\},\]
and
\[S_3=\#\{x_1,x_2,x_3,x_4\le q:\,q|x_1^2+x_2^2-x_3^2-x_4^2\}.\]
We shall relate $S_1$ and $S_2$ to $S_3$, using exponential sums.
If we write \linebreak $e_q(m):=\exp(2\pi im/q)$ we can use a standard
manipulation to show that
\[S_1=q^{-4}\sum_{\lb{b}} S(\b{b};q)T(b_1;N,q)T(b_2;N,q)T(b_3;N,q)T(b_4;N,q),\]
where $\b{b}=(b_1,b_2,b_3,b_4)$ runs over vectors modulo $q$,
and the sums $S(\b{b};q)$ and $T(b;N,q)$ are given by
\[S(\b{b};q)=\twosum{\lb{x}\mod{q}}{q|x_1^2+x_2^2-x_3^2-x_4^2}
e_q(\b{b}.\b{x})\]
and
\[T(b;N,q)=\sum_{x\le N}e_q(-bx)\ll\min(N\,,\,||b/q||^{-1}).\]
Similarly we find that
\[S_2=q^{-2}\sum_{b_1,b_2} S((b_1,b_2,0,0);q)T(b_1;N,q)T(b_2;N,q)\]
and
\[S_3=S((0,0,0,0);q).\]
Since $T(0;N,q)=N$ we see that the terms corresponding to 
$\b{b}=\b{0}$ cancel in (\ref{kk}),
and it remains to estimate the contribution to $S_1$ and $S_2$ arising
from terms with $\b{b}\not=\b{0}$.  We shall write
\[S^{(i)}=
q^{-4}\sum_{\lb{b}}|S(\b{b};q)T(b_1;N,q)T(b_2;N,q)T(b_3;N,q)T(b_4;N,q)|\]
where the sum is over vectors with $|b_j|\le q/2$, precisely $i$ of
which are non-zero.  Then
\begin{equation}\label{eno}
S\ll \sum_{i=1}^4 S^{(i)}.
\end{equation}

The sums $S(\b{b};q)$ satisfy a product rule
\[S(\b{b};q_1q_2)=S(\b{b};q_1)S(\b{b};q_2),\;\;\;(q_1,q_2)=1\]
and a trivial bound
\[|S(\b{b};q)|\le q^4.\]
Moreover when $q$ is an odd prime $p$ a standard evaluation shows that
\[S(\b{b};p)=p^3+p^2-p\] 
when $p|\b{b}$; that $S(\b{b};p)=p^2-p$ when 
$p|b_1^2+b_2^2-b_3^2-b_4^2$ but $p\nmid\b{b}$; and that
$S(\b{b};p)=-p$ if $p\nmid b_1^2+b_2^2-b_3^2-b_4^2$.
We may therefore decompose $q$ into coprime factors $q=q_1q_2q_3q_4$
such that $q_2,q_3,q_4$ are odd and square-free, with $q_2|\b{b}$ and 
$q_3|b_1^2+b_2^2-b_3^2-b_4^2$.  Moreover we will have
\[S(\b{b};q)\ll q_1^4q_2^{3+\eta}q_3^2q_4.\]
There are $O(q^{\eta})$ possible factorizations $q=q_1q_2q_3q_4$.
Thus
\begin{equation}\label{B}
S^{(i)}\ll q^{-4+\eta}\max_{q_1,q_2,q_3,q_4}q_1^4q_2^{3+\eta}q_3^2q_4
S_i(q_1,q_2,q_3,q_4),
\end{equation}
in which
\begin{eqnarray*}
S_i(q_1,q_2,q_3,q_4)&=&
\sum_{\lb{b}}|T(b_1;N,q)T(b_2;N,q)T(b_3;N,q)T(b_4;N,q)|\\
&\ll&\sum_{\lb{b}}\prod_{j=1}^4\min(N\,,\,q|b_j|^{-1}),
\end{eqnarray*}
where $\b{b}$ runs over vectors in the range $|b_j|\le q/2$, precisely
$i$ of which are non-zero, and for which
\[q_2|\b{b},\;\;\; q_3|b_1^2+b_2^2-b_3^2-b_4^2.\]

We shall discuss the case of $S^{(4)}$ in detail, the other sums
being treated similarly.  For any fixed choice of $\pm$ signs we write
\[K(C_1,C_2,C_3,C_4;q_3)=\#\{\b{c}\in\Z^4: 
q_3|c_1^2\pm c_2^2\pm c_3^2\pm c_4^2,\;
|c_j|\le C_j,\, (1\le j\le 4)\},\]
where we shall assume that $1\le C_j\le q$ for all $j$.
To estimate this we shall suppose that $C_1\le C_2\le C_3\le C_4$. 
If $c_1^2\pm c_2^2\pm c_3^2\pm c_4^2=q_3k$, say, then
$k\ll C_4^2/q_3$.  For each value of $k$ one sees that $c_1$ and $c_2$ 
determine $O(q^{\eta})$ pairs $c_3,c_4$, unless $c_1^2\pm c_2^2=q_3k$, 
in which case $c_1$ and $c_3$ determine $O(1)$ pairs $c_2,c_4$.  
Thus there are $O(C_1C_3q^{\eta})$ possibilities for each value of
$k$, so that
\[K(C_1,C_2,C_3,C_4;q_3)\ll (1+C_4^2/q_3)C_1C_3q^{\eta}.\]
For an alternative estimate we observe that the congruence
$c_4^2\equiv n\mod{q_3}$ 
has $O(q_3^{\eta})$ solutions modulo $q_3$, whence
\[K(C_1,C_2,C_3,C_4;q_3)\ll (1+C_4/q_3)C_1C_2C_3q^{\eta}.\]
To put these bounds into a more convenient form we write
$C:=C_1C_2C_3C_4$ and observe that by combining our estimates we have
\[K(C_1,C_2,C_3,C_4;q_3)\ll
\{C_1C_3+q_3^{-1}C+\min(q_3^{-1}C_1C_3C_4^2\,,\, C_1C_2C_3)\}
q^{\eta}.\]
Since $C_1\le C_2\le C_3\le C_4$ we have
$C_1C_3\le\sqrt{C}$ and 
\begin{eqnarray*}
\min(q_3^{-1}C_1C_3C_4^2\,,\, C_1C_2C_3)&\le& 
(q_3^{-1}C_1C_3C_4^2)^{1/2}(C_1C_2C_3)^{1/2}\\
&=&q_3^{-1/2}C_1C_2^{1/2}C_3C_4\\
&\le& q_3^{-1/2}CC_0^{-1/2},
\end{eqnarray*}
where we have written $C_0=\min C_i$.  It follows that
\begin{eqnarray*}
K(C_1,C_2,C_3,C_4;q_3)&\ll&
\{C^{1/2}+q_3^{-1}C+q_3^{-1/2}CC_0^{-1/2}\}q^{\eta}\\
&\ll& C\{C_0^{-2}+q_3^{-1}+q_3^{-1/2} C_0^{-1/2}\}q^{\eta}.
\end{eqnarray*}

We can now bound $S_4(q_1,q_2,q_3,q_4)$. We write $b_j=q_2c_j$ and
decompose the range for each
$c_j$ into intervals either of the shape $|c_j|\le C_j=qq_2^{-1}N^{-1}$, or of
the form $C_j/2\le|c_j|\le C_j$ with $C_j\ge qq_2^{-1}N^{-1}$.  
Ranges with $C_j<1$ will not arise, since all $b_j$ are non-zero for $S^{(4)}$.
There will be $\ll\log^4 q\ll q^{\eta}$
sets of ranges in total, on each of which we will have
\[\prod_{j=1}^4\min(N\,,\,q|b_j|^{-1})\ll q^4q_2^{-4}C^{-1}.\]
Moreover, since $C_0\ge qq_2^{-1}N^{-1}\ge q^{1/3}q_2^{-1}$, we find that
\begin{eqnarray*}
K(C_1,C_2,C_3,C_4;q_3)&\ll&
C\{q^{-2/3}q_2^{2}+q_3^{-1}+q_3^{-1/2}q^{-1/6}q_2^{1/2}\}q^{\eta}\\
&\ll& Cq_2^{4/3}q_3^{-2/3}q^{\eta}.
\end{eqnarray*}
It now follows that
\begin{eqnarray*}
S_4(q_1,q_2,q_3,q_4)&\ll& 
q^{4+\eta}q_2^{-4}C^{-1}.Cq_2^{4/3}q_3^{-2/3}q^{\eta}\\
&=& q^{4+2\eta}q_2^{-8/3}q_3^{-2/3}
\end{eqnarray*}
whence (\ref{B}) yields
\begin{eqnarray*}
S^{(4)}&\ll& 
q^{-4+\eta}q_1^4q_2^{3+\eta}q_3^2q_4.q^{4+2\eta}q_2^{-8/3}q_3^{-2/3}\\
&\ll& q^{4/3+4\eta}q_1^{3}.
\end{eqnarray*}

Similar arguments show that
\begin{eqnarray*}
\lefteqn{\#\{\b{c}\in\Z^3: q_3|c_1^2\pm c_2^2\pm c_3^2,\;
|c_j|\le C_j,\, (1\le j\le 3)\}}\hspace{6cm}\\
&\ll& C_1C_2(1+q_3^{-1}C_3)q_3^{\eta}\\
&\ll& (C_1C_2C_3)q_2^{2/3}q_3^{-1/3}q^{\eta},
\end{eqnarray*}
whence
\[S^{(3)}\ll Nq^{2/3+4\eta}q_1^3\ll q^{4/3+4\eta}q_1^3;\]
that
\begin{eqnarray*}
\#\{\b{c}\in\Z^2: q_3|c_1^2\pm c_2^2,\;
|c_j|\le C_j,\, (1\le j\le 2)\}&\ll& C_1(1+q_3^{-1}C_2)q_3^{\eta}\\
&\ll& (C_1C_2)q_2^{2/3}q_3^{-1/3}q^{\eta},
\end{eqnarray*}
whence
\[S^{(2)}\ll N^2q^{-1/3+4\eta}q_1^3\ll q^{1+4\eta}q_1^3;\]
and that
\[\#\{c\in\Z: q_3|c^2,\; |c|\le C\}\ll C/q_3\]
whence
\[S^{(1)}\ll N^3q^{-1+4\eta}q_1^3\ll q^{1+4\eta}q_1^3.\]
In view of (\ref{eno}) these estimates suffice for the proof of
Lemma \ref{dispa}.
\bigskip

We proceed to deduce Lemma \ref{disp} from Lemma \ref{dispa}.  
If $0\le M\le N$ then
\begin{eqnarray*}
\lefteqn{A(N+M,q,c)-A(N,q,c)}\hspace{1cm}\\
&\ll& \#\{u\in(N,N+M],\,x\le N+M: u^2-x^2\equiv\pm c\mod{q}\},
\end{eqnarray*}
whence
\begin{eqnarray*}
\lefteqn{\sum_{c=1}^q \{A(N+M,q,c)-A(N,q,c)\}^2}\hspace{1cm}\\
&\ll&\#\{u,v\in(N,N+M],\,x,y\le N+M:\, q|u^2-x^2\pm(v^2-y^2)\}.
\end{eqnarray*}
If $u^2-x^2\pm(v^2-y^2)=kq\ll N^2$ then each set of values $u,v,k$
determines $O(q^{\eta})$ pairs $x,y$, except when the $\pm$ sign is
negative and $u^2-v^2=kq$. Thus
there are $O(M^2N^2q^{-1+\eta})$ solutions $u,v,x,y$ with 
$u^2-v^2\not=kq$.  When
$u^2-v^2=kq$ we have $k\ll MNq^{-1}$.  Thus there are
$O(MNq^{-1+\eta})$ pairs $u,v$ corresponding to non-zero values of
$k$, and $M$ pairs for $k=0$.  To each such pair $u,v$ with
$u^2-v^2=kq$ there correspond $O(N)$ pairs $m=n$.  We therefore 
obtain the bound
\[\sum_{c=1}^q \{A(N+M,q,c)-A(N,q,c)\}^2\ll
M^2N^2q^{-1+\eta}+MN^2q^{-1+\eta}+MN.\]
In particular, taking $M=N=q$, we have
\[\sum_{c=1}^qA_0(q,c)^2\ll q^{3+\eta}.\]
It follows that
\begin{eqnarray*}
\lefteqn{\sum_{c=1}^q \max_{0\le H\le M}|\Delta(N+H,q,c)-\Delta(N,q,c)|^2}
\hspace{1cm}\\
&\ll & \sum_{c=1}^q \max_{0\le H\le M}\{A(N+H,q,c)-A(N,q,c)\}^2\\
&&\hspace{2cm}\mbox{}
+\max_{0\le H\le M}\frac{((N+H)^2-N^2)^2}{q^4}\sum_{c=1}^q A_0(q,c)^2\\
&=&\sum_{c=1}^q\{A(N+M,q,c)-A(N,q,c)\}^2\\
&&\hspace{2cm}\mbox{}
+\frac{((N+M)^2-N^2)^2}{q^4}\sum_{c=1}^q A_0(q,c)^2\\
&\ll & q^{\eta}\{M^2N^2q^{-1}+MN\}\\
&\ll& q^{\eta}\{M^2q^{1/3}+Mq^{2/3}\},
\end{eqnarray*}
if $N\le q^{2/3}$.

To handle
\[\Delta^*(q,c)=\max_{U\le q^{2/3}}\Delta(U,q,c)\]
we cover the available range for $U$ with $O(q^{2/3}M^{-1})$
sub-intervals $N_j\le U\le N_j+M$, whence Lemma \ref{dispa} yields
\begin{eqnarray*}
\lefteqn{\sum_{c=1}^q\Delta^*(q,c)^2}\\
&\ll& \sum_{N_j}\sum_{c=1}^q\{\Delta(N_j,q,c)^2+\max_{0\le H\le M}
|\Delta(N_j+H,q,c)-\Delta(N_j,q,c)|^2\}\\
&\ll_{\eta}& q^{2/3}M^{-1}q^{4\eta}\{q^{4/3}q_1^3+q^{1/3}M^2+q^{2/3}M\}.
\end{eqnarray*}
The choice $M=q^{1/2}$ then results in the estimate required for Lemma
\ref{disp}.

\section{Proof of Theorem \ref{main}}

On writing $n-m=u,\, n+m=v$ we find that
\[R_{\alpha}(N,X)=N^{-1}\#\{m<n\le N:||\alpha(m^2-n^2)||\le
XN^{-1}\}\]
\[=N^{-1}\sum_{u\le N}\#\{v\equiv u\mod{2}:\, u<v\le 2N-u,\,
||\alpha uv||\le XN^{-1}\}.\]
When $u$ is even we write $2u$ in place of $u$ and put $v=2x$ to find that the
corresponding contribution is
\begin{eqnarray}\label{f}
\lefteqn{N^{-1}\sum_{u\le N/2}\#\{x\in \N:\, u<x\le N-u,\,
||4\alpha ux||\le XN^{-1}\}}\hspace{1cm}\nonumber\\\
&=&
N^{-1}\sum_{u\le N/2}\{R(N-u,4\alpha u,X/N)-R(u,4\alpha u,X/N)\},
\end{eqnarray}
say, where
\[R(M,\beta,\delta):=\#\{x\in \N:\, x\le M,\,
||\beta x||\le\delta\}.\]
To handle odd values $u$ we count integers $v\equiv 1\pmod{2}$ by first
considering the contribution from all $v$, and then subtracting the
contribution from even $v$.  This leads to a total
\begin{eqnarray*}
\lefteqn{N^{-1}\sum_{u\le N,\,2\nmid u}
\{R(2N-u,\alpha u,X/N)-R(u,\alpha u,X/N)\}}\hspace{1cm}\\
&&-N^{-1}\sum_{u\le N,\, 2\nmid u}
\{R(N-u/2,2\alpha u,X/N)-R(u/2,2\alpha u,X/N)\}.
\end{eqnarray*}

We proceed to estimate $R(M,\beta,\delta)$.
Let $\b{u},\b{v}\in\R^2$ be the vectors 
\[\b{u}=\left(\sqrt{\frac{\delta}{M}}\, ,
  \,\beta\sqrt{\frac{M}{\delta}}\right), \;\;\;\mbox{and}\;\;\;
\b{v}=\left(0\, , \, -\sqrt{\frac{M}{\delta}}\right).\]
Then for $\delta\in(0,1)$ we have
\begin{eqnarray*}
1+2R(M,\beta,\delta)&=&
\#\{(x,y)\in\Z^2:\, |x|\le M,\, |\beta x-y|\le\delta\}\\
&=&\#\{(x,y)\in\Z^2:\, x\b{u}+y\b{v}\in
[-\sqrt{M\delta}\,,\,\sqrt{M\delta}\,]^2\}.
\end{eqnarray*}
The vectors $\b{u},\b{v}$ generate a lattice of determinant 1.  Hence 
\begin{equation}\label{lat0}
\#\{(x,y)\in\Z^2:\, x\b{u}+y\b{v}\in [\, -S,S\, ]^2\}
=(2S)^2+O(S/\lambda_1)+O(1), 
\end{equation}
where $\lambda_1$ is the first successive
minimum of the lattice, that is to say the length of the shortest
non-zero vector in the lattice.  In our case we find that
\[1+2R(M,\beta,\delta)=4M\delta+O(M^{1/2}\delta^{1/2}\lambda_1^{-1})
+O(1).\]
We have 6 different pairs $(M,\beta)=(N-u,4\alpha
u),\ldots,(u/2,2\alpha u)$ to consider, each with a corresponding
value for $\lambda_1$.  We write $\lambda_0$ for the smallest of these
6 values, and split the available range for $M^{1/2}\delta^{1/2}\lambda_0^{-1}$
into dyadic intervals
\[4E< M^{1/2}\delta^{1/2}\lambda_0^{-1}\le 8E.\]
For values $u$ for which $E\le X^{1/2}$ the total 
contribution of the error terms to $R_{\alpha}(N,X)$ is clearly $O(X^{1/2})$.
The choice of $X^{1/2}$ as the point at which we split the range for
$E$ is not optimal, but is adequate for our purposes.

In the remaining case $E>X^{1/2}$ there will be coprime
integers $x,y$ for which $|x|\le M/(4E)$ and $|\beta
x-y|\le \delta/(4E)$.  It follows that we will have
\begin{eqnarray*}
R(N-u,4\alpha u,X/N)&=&2X(1-u/N)+O(E),\\
R(u,4\alpha u,X/N)&=&2Xu/N+O(E),\\
R(2N-u,\alpha u,X/N)&=&2X(2-u/N)+O(E),\\
R(u,\alpha u,X/N)&=&2Xu/N+O(E)\\
R(N-u/2,2\alpha u,X/N)&=&X(2-u/N)+O(E),\\
R(u/2,2\alpha u,X/N)&=&Xu/N+O(E),
\end{eqnarray*}
and that there is a coprime pair $x,y$ satisfying one of
\[|x|\le \frac{N}{4E},\;\;\; |4\alpha ux -y|\le \frac{X}{4NE}\]
or
\[|x|\le \frac{N}{2E},\;\;\; |\alpha ux -y|\le \frac{X}{4NE}\]
or
\[|x|\le \frac{N}{4E},\;\;\; |2\alpha ux -y|\le \frac{X}{4NE}\]
respectively.  To simplify matters we replace $x$ by $x'=4x$ in the first 
case and by $x'=2x$ in the third, and then remove a factor $(x',y)=2$
or 4 if necessary.  We deduce in each case that there is a coprime
pair with $|x|\le N/E$ and
$|\alpha ux-y|\le X/(NE)$.  Since they are coprime, $x$ and $y$ cannot both
vanish.  Indeed, since $X\le \log N$ we will have $X/(NE)<1$ whence it
is clear that $x$ cannot vanish.  It 
follows that the total contribution of the error terms $O(E)$ to
$R_{\alpha}(N,X)$, arising from an individual value of $E\ge X^{1/2}$, is
\[\ll \frac{E}{N}
\#\{(u,x,y)\in\N^2\times\Z:\, u\le N,\,x\le \frac{N}{E},\, (x,y)=1,\,
|\alpha ux-y|\le \frac{X}{NE}\}.\]

We now calculate the contribution from the main term of
$R(M,\beta,\delta)$.  For (\ref{f}) this is
\[\frac{1}{N}\sum_{u\le N/2}\{2X(1-u/N)-2Xu/N\}=
\frac{4X}{N^2}\sum_{u\le N/2}(N/2-u)=\frac{X}{2}+O(\frac{X}{N}).\]
Similarly the odd values of $u$ contribute $X/2+O(X/N)$.
We therefore conclude as follows.
\begin{lemma}\label{lat}
If $1\le X\le\log N$ then
\[R_{\alpha}(N,X)=X+O(X^{1/2})+O\left(N^{-1}\sum_{X^{1/2}\le E=2^k\le N}E
V^*(N,\frac{N}{E};\frac{X}{NE})\right)\]
where
\begin{eqnarray*}
\lefteqn{V^{*}(A,B;\Delta):=}\hspace{1cm}\\
&&\#\{(u,x,y)\in\N^2\times\Z:\, u\le A,\,x\le B,\, 
(x,y)=1,\, |\alpha ux-y|\le \Delta\}.
\end{eqnarray*}
\end{lemma}
It is already clear here that our approach cannot provide an
asymptotic evaluation for $R_{\alpha}(N,X)$ unless
$X\rightarrow\infty$.  The error term $O(1)$ in (\ref{lat0}) will
produce at least a corresponding error $O(1)$ for $R_{\alpha}(N,X)$.  
Any sharper estimate
would appear to require information on the way the shape of our
lattice varies with the parameter $u$.

From now on we shall focus on the second error term above.  We write
$(u,y)=f$ and suppose that $f$ lies in a dyadic range $F\le f<2F$.
Given such an $f$, if $u=fu_0$ and $y=fy_0$ then $u_0\le N/F$ and
$|\alpha u_0x-y_0|\le X/(NEF)$.
Moreover each pair $u_0,y_0$ can correspond to at most $F$ pairs
$u,y$, since we are assuming that $F\le f<2F$.  Our error term is
therefore
\begin{equation}\label{err}
\ll N^{-1}\sum_{X^{1/2}\le E=2^k\le N}\;\sum_{1\le F=2^h\le N}EF
V(\frac{N}{E},\frac{N}{F};\frac{X}{NEF})
\end{equation}
where
\begin{eqnarray*}
\lefteqn{V(A,B;\Delta):=}\hspace{1cm}\\
&&\#\{(a,b,z)\in\N^2\times\Z:\, a\le A,\,b\le B,\, (ab,z)=1,\, 
|\alpha ab-z|\le \Delta\}.
\end{eqnarray*}

Our strategy for tackling $V(A,B;\Delta)$ is based on the following
lemma.
\begin{lemma}\label{sf}
Let $E$ be an ellipse centred at the origin, of area $A(E)$.  Then the
number of coprime integer pairs $(x,y)\in E$ is $O(1+A(E))$.
\end{lemma}
This easy result may be found in the author's work \cite[Lemma 2]{sf},
for example.

If we fix $a$, say, then the ellipse
\[\{(b,z)\in\R^2: B^{-2}b^2+\Delta^{-2}|\alpha ab-z|^2\le 2\}\]
has area $\ll B\Delta$, and we deduce that $V(A,B;\Delta)\ll
AB\Delta+A$.  By symmetry we then have the bound  
\begin{equation}\label{simp}
V(A,B;\Delta)\ll AB\Delta+\min(A,B).  
\end{equation}
In our application the contribution from the term
$AB\Delta$ is usually satisfactory, but the effect of the second term is
likely to be too large unless $A=N/E$ and $B=N/F$ have very different
sizes.  To circumvent this difficulty we shall use a delicate arithmetic
trick, which is the key to our attack on Theorem \ref{main}.  

It will be convenient to assume that $A\le B$, as we may, by
symmetry.  We take
parameters $P_1\ge P_0\ge 1$ and consider prime factors $p$ of $ab$ in
the range $P_0<p\le P_1$.  Thus we will need to consider separately
\[V_1:=\#\{(a,b,z)\in\N^2\times\Z:\, a\le A,\,b\le B,\, (ab,z\Pi)=1,\, 
|\alpha ab-z|\le \Delta\},\]
where
\[\Pi:=\prod_{P_0<p\le P_1}p.\]
Here Lemma \ref{sf} shows that
\begin{eqnarray*}
V_1&\le&\sum_{a\le A, (a,\Pi)=1}\#\{(b,z)\in\Z^2:\, b\le B,\, (b,z)=1,\, 
|\alpha ab-z|\le \Delta\}\\
&\ll& (B\Delta+1)\#\{a\le A:\, (a,\Pi)=1\}.
\end{eqnarray*}
The number of available integers $a$ may be estimated using a standard
sieve bound.  According to Theorem 2.2 of Halberstam and Richert \cite{HR},
for example, one has
\[\#\{a\le A:\, (a,\Pi)=1\}\ll  A\prod_{P_0< p\le P_1}(1-1/p)
\ll  A\frac{\log P_0}{\log P_1}\]
providing that $P_1\le A$.  This yields 
the following lemma.
\begin{lemma}\label{v1}
If $1\le P_0\le P_1\le A$ we have
\[V_1\ll AB\Delta+A\frac{\log P_0}{\log P_1}.\]
\end{lemma}

If $ab$ does have a prime factor $p$ in the range $P_0<p\le P_1$ we
may choose the smallest such prime $p$, and classify the corresponding
triples $a,b,z$ according to the dyadic range $P=2^k<p\le 2P$ in which $p$
lies.  We write $V_2(P)$ for the corresponding contribution to
$V(A,B;\Delta)$ and write
\[V_2:=\sum_{P_0\le P=2^k\le P_1}V_2(P).\]
If $p|a$ we set $a'=a/p$ and $b'=bp$, while if $p\nmid a$ we
will have $p|b$, and we set $a'=ap$ and $b'=b/p$.  It follows that
$|\alpha a'b'-z|\le \Delta$, and that either $a'\le A/P,\, b'\le 2BP$
or $a'\le 2AP,\, b'\le B/P$.  Moreover, if we are given a triple
$a',b',z$ counted by $V(A/P,2BP;\Delta)$ then it determines the prime
$p$, which will be the smallest prime $p>P_0$ dividing $a'b'=ab$.
Knowing $p$ one may then find the pair $a,b$ which must either be
$a=a'p,\, b=b'/p$ or $a=a'/p,\, b=b'p$.  It follows that each triple
$a',b',z$ counted by 
\[V(A/P,2BP;\Delta)+V(2AP,B/P;\Delta)\]
arises from
at most 2 triples $a,b,z$ counted by $V_2(P)$.  We may now use
(\ref{simp}) to deduce that
\[V_2(P)\ll AB\Delta+\min(A/P,BP)+\min(AP,B/P)\ll AB\Delta+B/P,\]
from which we obtain the following lemma.
\begin{lemma}\label{v2}
If $1\le P_0\le P_1$ we have
\[V_2\ll AB\Delta(\log P_1)+P_0^{-1}B.\]
\end{lemma}

We now use Lemmas \ref{v1} and \ref{v2} to estimate the contribution
to (\ref{err}) from terms with $EF\le (\log N)^{5/4}$.  We
choose 
\[P_0=E^2F^2,\;\;\; P_1=6\exp((EF)^{3/4}),\] 
so that $1\le P_0\le P_1\le A$ for large enough $N$.
The terms $AB\Delta$ and $AB\Delta(\log P_1)$ in Lemmas \ref{v1} and
\ref{v2} then contribute
\begin{eqnarray*}
&\ll &N^{-1}\sum_{X^{1/2}\le E=2^k\le N}\;\sum_{1\le F=2^h\le N}
EF\frac{N^2}{EF}\frac{X}{NEF}(EF)^{3/4}\\
&\ll& \sum_{X^{1/2}\le E=2^k\le N}\;\sum_{1\le F=2^h\le N}
\frac{X}{(EF)^{1/4}}\\
&\ll& X^{7/8}.
\end{eqnarray*}
Moreover the error term $A(\log P_0)/(\log P_1)$
in Lemma \ref{v1} produces
\[\ll N^{-1}\sum_{X^{1/2}\le E=2^k\le N}\;\sum_{1\le F=2^h\le N}
EF\min(\frac{N}{E}\,,\,\frac{N}{F})\frac{\log(EF)}{(EF)^{3/4}}\ll 1,\]
since $\min(N/E,N/F)\le N(EF)^{-1/2}$.
Finally, the error term $P_0^{-1}B$ occurring in Lemma~\ref{v2} produces
\[\ll N^{-1}\sum_{X^{1/2}\le E=2^k\le N}\;\sum_{1\le F=2^h\le N}
EF\max(\frac{N}{E}\,,\,\frac{N}{F})E^{-2}F^{-2}\ll 1.\]
Thus those terms with $EF\le (\log N)^{5/4}$ make a satisfactory
contribution in Theorem \ref{main}.

Up to this point we have made no use of the Diophantine approximation
properties of $\alpha$, but it is time to bring these into play.  In
order to clarify the rationale behind our choice of the various
exponents which will occur, we introduce constants
$\beta,\gamma\in(0,1)$ on which we will impose certain constraints as the
argument progresses, and which will eventually be specified in
(\ref{ass1}).  To begin with we assume that $\alpha$ satisfies
\[\left|\alpha-\frac{a}{q}\right|\ge \frac{1}{\kappa q^{2+\beta}}\]
for every fraction $a/q$.  In particular, if  
$V(N/E,N/F;X/(NEF))$ counts $x,y,z$, so that 
$|\alpha-z/xy|\le X/(NEFxy)$, we deduce
that
\[(xy)^{1+\beta}\ge NEF/(\kappa X).\]
Hence $V(N/E,N/F;X/(NEF))=0$ unless 
\[(N^2/EF)^{1+\beta}\ge NEF/(\kappa X).\]
We therefore assume from
now on that
\begin{equation}\label{EF}
(\log N)^{5/4}\le EF\le 
(\kappa X)^{1/(2+\beta)}N^{(1+2\beta)/(2+\beta)}.
\end{equation}

Let
\begin{equation}\label{Qc}
Q=\left(\frac{N^2}{EF}\right)^{1-\gamma}
\end{equation}
and apply Dirichlet's Approximation Theorem to obtain coprime integers $a,q$
with 
\[\left|\alpha-\frac{a}{q}\right|\le \frac{1}{qQ},\;\;\; 1\le q\le Q.\]
It follows of course that 
\begin{equation}\label{qb}
(\kappa^{-1}Q)^{1/(1+\beta)}\le q\le Q.
\end{equation}
Now if $||\alpha xy||\le X/(NEF)$ with $x\le N/E$ and $y\le N/F$ then
\[||axy/q||\le X/(NEF)+N^2/(EFqQ).\]
whence $axy\equiv r\mod{q}$ for
some integer $r$ with 
\[|r|\le qX/(NEF)+N^2/(EFQ).\]
It follows that
\begin{equation}\label{e4}
V(\frac{N}{E},\frac{N}{F};\frac{X}{NEF})\le\sum_{|r|\le qX/(NEF)+N^2/(EFQ)}
\;\twosum{n\le N^2/(EF)}{an\equiv r\mod{q}}d(n),
\end{equation}
where $d(n)$ is the divisor function.
The reader should observe that there is a loss at this point, in replacing
\[\#\{x,y:\,x\le N/E,\, y\le N/F,\, xy=n\}\]
by $d(n)$.  This loss is of order $\log N$, and is only acceptable
since we are now in the case in which $EF$ is of larger order than $\log N$.

For the case in which $(s,q)=1$ it was shown by Linnik and Vinogradov
\cite{LV} that
\[\twosum{m\le M}{m\equiv
  s\mod{q}}d(m)\ll_{\gamma}\frac{\phi(q)}{q^2}M\log M,\]
providing that $q\le M^{1-\gamma}$ for some constant $\gamma>0$.  In
general, if $(s,q)=h$ say, one may write $s=hs',\, q=hq'$ and $m=hm'$,
so that $q'\le (M/h)^{1-\gamma}$ and $d(m)\le d(h)d(m')$.  Then
\begin{eqnarray*}
\twosum{m\le M}{m\equiv s\mod{q}}d(m)&\le & 
d(h)\twosum{m'\le M/h}{m'\equiv s'\mod{q'}}d(m')\\
&\ll_{\gamma}& d(h)\frac{\phi(q')}{{q'}^2}\frac{M}{h}\log M\\
&\ll_{\gamma}& d(h)q^{-1}M\log M.
\end{eqnarray*}
If we sum for $|s|\le S$ we find that
\begin{eqnarray*}
\sum_{|s|\le S}d((s,q))&\le&\sum_{k|q}d(k)\#\{s:\,|s|\le S,\, k|s\}\\
&\ll&\sum_{k|q}d(k)(1+S/k)\\
&\ll& d^2(q)+S\prod_{p|q}(1+2/p+O(p^{-2}))\\
&\ll& d^2(q)+S\frac{\sigma(q)^2}{q^2}\\
&\ll& d^2(q)+S(\log\log q)^2.
\end{eqnarray*}
It follows from (\ref{e4}) that
\begin{eqnarray}\label{vb}
\lefteqn{V(\frac{N}{E},\frac{N}{F};\frac{X}{NEF})}\nonumber\\
&\ll& \{d^2(q)+(\frac{qX}{NEF}+\frac{N^2}{EFQ})(\log\log q)^2\}
\frac{N^2}{qEF}\log N.
\end{eqnarray}

We begin by examining the first term $d^2(q)N^2(qEF)^{-1}\log N$.  In
view of (\ref{qb}) and (\ref{Qc}) this is at most
\[\kappa^{1/(1+\beta)}d^2(q)(\log N)\frac{N^2}{EF}Q^{-1/(1+\beta)}
=\kappa^{1/(1+\beta)}d^2(q)(\log N)
\left(\frac{N^2}{EF}\right)^{(\beta+\gamma)/(1+\beta)}.\]
When we multiply by $N^{-1}EF$ and sum over dyadic ranges subject to
(\ref{EF}) we see that the contribution to (\ref{err}) is
\[\ll
\kappa^{1/(1+\beta)}
(\kappa X)^{(1-\gamma)/(1+\beta)(2+\beta)}
d^2(q)(\log N)^2N^{-\phi_1}\]
with
\[\phi_1=
1-(1+\beta)^{-1}\{2(\beta+\gamma)+(1-\gamma)\frac{1+2\beta}{2+\beta}\}.\]

Turning to the second term on the right of (\ref{vb}), we see that the
overall contribution to the error terms in Theorem \ref{main} is
\[X(\log\log q)^2(\log N)\sum_{EF\ge(\log N)^{5/4}}(EF)^{-1}\ll
X(\log\log N)^3(\log N)^{-1/4},\]
which is $O(X^{7/8})$. This is satisfactory for the theorem.

Finally, the third term on the right of (\ref{vb}) is
\begin{eqnarray*}
&\ll &(\log\log N)^2(\log N)\left(\frac{N^2}{EF}\right)^2
\kappa^{1/(1+\beta)}Q^{-(2+\beta)/(1+\beta)}\\
&=&\kappa^{1/(1+\beta)}(\log\log N)^2(\log N)
\left(\frac{N^2}{EF}\right)^{(2\gamma+\beta\gamma+\beta)/(1+\beta)}.
\end{eqnarray*}
We impose the condition that $\beta,\gamma<1/3$, which ensures that the
exponent $(2\gamma+\beta\gamma+\beta)/(1+\beta)$ is less than 1.
Now, when we multiply by $N^{-1}EF$ and sum over dyadic ranges 
subject to (\ref{EF}), we get an overall contribution
\[\ll
\kappa^{1/(1+\beta)}(\kappa X)^{(1-2\gamma-\beta\gamma)/(1+\beta)(2+\beta)}
(\log\log N)^2(\log N)^2N^{-\phi_2},\]
with
\[\phi_2=1-(1+\beta)^{-1}\{2\beta+3\gamma+\frac{1+2\beta}{2+\beta}\}.\]

If $\gamma$ were equal to zero we would have
\[\phi_1=\phi_2=\frac{1-3\beta-\beta^2}{(1+\beta)(2+\beta)}.\]  
We will first choose $\beta$ so as to make this value positive, and 
then select a sufficiently small $\gamma$ so that
$\phi_1,\phi_2>0$. With this in mind we specify
\begin{equation}\label{ass1}
\beta=1/4 ,\;\;\; \gamma=1/40,
\end{equation}
from which the assertion of Theorem \ref{main} follows.

\bigskip
\bigskip

Mathematical Institute,

24--29, St. Giles',

Oxford

OX1 3LB

UK
\bigskip

{\tt rhb@maths.ox.ac.uk}

\end{document}